\documentclass[a4paper,11pt]{amsart}
\usepackage{graphicx}
 \usepackage{mathptmx}
\usepackage{amsmath}
\usepackage{amssymb}
\usepackage{enumitem}
\usepackage{xcolor}

\newmuskip\pFqmuskip

\newcommand*\pFq[6][8]{%
  \begingroup 
  \pFqmuskip=#1mu\relax
  \mathcode`=\string"8000
  \begingroup\lccode`\~=`\,
  \lowercase{\endgroup\let~}\pFqcomma
  F^{#2}_{#3}{\left(\genfrac..{0pt}{}{#4}{#5}\bigg|#6\right)}%
  \endgroup
}
\newcommand{\pFqcomma}{\mskip\pFqmuskip}

\newtheorem{theorem}{Theorem}[section]

\newtheorem{corollary}[theorem]{Corollary}

\newtheorem{remark}[theorem]{Remark}

\begin{document}

\title[Infinite series whose terms involve truncated degenerate exponentials]{A note on infinite series whose terms involve truncated degenerate exponentials}

\author{Dae San  Kim }
\address{Department of Mathematics, Sogang University, Seoul 121-742, Republic of Korea}
\email{dskim@sogang.ac.kr}
\author{Hyekyung Kim}
\address{Department Of Mathematics Education, Daegu Catholic University, Gyeongsan 38430, Republic of Korea}
\email{hkkim@cu.ac.kr}
\author{Taekyun  Kim}
\address{Department of Mathematics, Kwangwoon University, Seoul 139-701, Republic of Korea}
\email{tkkim@kw.ac.kr}


\maketitle

\begin{abstract}
The degenerate exponentials play an important role in recent study on degenerate versions of many special numbers and polynomials, the degenerate gamma function, the degenerate umbral calculus and the degenerate $q$-umbral calculus. The aim of this note is to consider infinite series whose terms involve truncated degenerate exponentials together with several special numbers and to find either their values or some other expressions of them as finite sums.
\end{abstract}

\section{Introduction} 

The degenerate exponentials play an important role in recent investigations on degenerate versions of many special numbers of polynomials (see [1]). Many of them are introduced by replacing the ordinary exponentials by the degenerate exponentials in their generating functions. These include the degenerate Stirling numbers of the second, the degenerate Bernoulli polynomials, the degenerate Euler polynomials, the partially degenerate Bell polynomials, and the degenerate central factorial numbers, and so on.
Not only that, the degenerate gamma function is introduced by replacing the ordinary exponential by the degenerate exponential in the integral representation of the usual gamma function (see [18]). Furthermore, as a degenerate version of the `classical' umbral calculus, the  $\lambda$-umbral calculus (also called degenerate umbral calculus) is developed again by making the same replacement in the generating function of the Sheffer sequences. As it turns  out, the degenerate umbral calculus (see [11]) is more convenient than the umbral calculus when dealing with degenerate special numbers and polynomials. In the same vein, the $\lambda$-$q$-umbral calculus (also called degenerate $q$-umbral calculus) is recently introduced by replacing the $q$-exponential by the $\lambda$-$q$-exponential (see [19]). In conclusion, we may say that study of degenerate versions has been very fruitful (see [2-4,10-19,21]).\par
The aim of this note is to consider several infinite series whose terms involve the truncated degenerate exponentials,\,\, $e_{\lambda}(y)-\frac{(1)_{k,\lambda}}{1}y-\cdots-\frac{(1)_{n,\lambda}}{n!}y^{n},\ (n\ge 0)$,\,\, and to find either their values or some other expressions of them as finite sums. Some of these infinite series also involve other special numbers, namely binomial coefficients, the generalized falling factorials (see \eqref{1-1}) and the degenerate Stirling numbers of the second kind (see \eqref{3}, \eqref{4}). \par

\vspace{0.1in}

For any $\lambda\in\mathbb{R}$, the degenerate exponentials are defined 
\begin{equation}
e_{\lambda}^{x}(t)=\sum_{n=0}^{\infty}\frac{(x)_{n,\lambda}}{n!}t^{n},\quad \mathrm{and}\quad e_{\lambda}(t)=e_{\lambda}^{1}(t)=\sum_{n=0}^{\infty}\frac{(1)_{n,\lambda}}{n!}t^{n},\label{1}	
\end{equation}
where the generalized falling factorials are given by
\begin{equation}
(x)_{0,\lambda}=1,\quad (x)_{n,\lambda}=x(x-\lambda)(x-2\lambda)\cdots(x-(n-1)\lambda),\quad (n\ge 1),\quad (\mathrm{see}\ [10,15]). \label{1-1}
\end{equation}
From \eqref{1}, we note that $\displaystyle\lim_{\lambda\rightarrow 0}e_{\lambda}^{x}(t)=e^{x}.\displaystyle$ The Stirling numbers of the second kind are given by 
\begin{equation}
x^{n}=\sum_{k=0}^{n}S_{2}(n,k)(x)_{k},\quad (n\ge 0),\quad (\mathrm{see}\ [10]),\label{2}	
\end{equation}
where $(x)_{k}=x(x-1)\cdots(x-k+1),\quad (k\ge 1),\quad (x)_{0}=0$. \par 
In [10], the degenerate Stirling numbers of the second kind are defined by 
\begin{equation}
	(x)_{n,\lambda}=\sum_{k=0}^{n}S_{2,\lambda}(n,k)(x)_{k},\quad (n\ge 0). \label{3}
\end{equation}
From \eqref{3}, we note that $\displaystyle\lim_{\lambda\rightarrow 0}S_{2,\lambda}(n,k)=S_{2}(n,k)\displaystyle$. \par 
By \eqref{3}, we easily get 
\begin{equation}
\frac{1}{k!}\Big(e_{\lambda}(t)-1\Big)^{k}=\sum_{n=k}^{\infty}S_{2,\lambda}(n,k)\frac{t^{n}}{n!},\quad (k\ge 0),\quad (\mathrm{see}\ [10,16,17]).\label{4}
\end{equation}
The backward difference operator $\bigtriangledown$ is defined as
\begin{equation}
\bigtriangledown f(x)=f(x)-f(x-1),\quad (\mathrm{see}\ [17]). \label{5}	
\end{equation}
 From \eqref{5}, we note that 
 \begin{equation}
 \binom{x-1}{n-1}= \bigtriangledown\binom{x}{n}=\binom{x}{n}-\binom{x-1}{n},\quad (n\ge 1). \label{6}
 \end{equation}
Thus, by \eqref{6}, we get 
\begin{equation}
\binom{x}{n}=\binom{x+1}{n}-\binom{x}{n-1},\quad (n\ge 0),\quad (\mathrm{see}\ [3,4,5,6,7]).\label{7}
\end{equation}
In addition, the degenerate Bell polynomials are defined by 
\begin{equation*}
\phi_{n,\lambda}(x)=e^{-x}\sum_{k=0}^{\infty}\frac{(k)_{n,\lambda}}{k!}x^{k}=\sum_{k=0}^{n}S_{2,\lambda}(n,k)x^{k},\quad (n\ge 0),\quad (\mathrm{see}\ [10,16]).
\end{equation*}

\section{Infinite series whose terms involve truncated degenerate exponentials}
In the section, we will consider infinite series whose terms involve truncated degenerate exponentials. We first observe that

\begin{align}
\frac{1}{x-1}&\bigg(e_{\lambda}(xy)-e_{\lambda}(y)\bigg)=\sum_{k=1}^{\infty}\frac{(1)_{k,\lambda}}{k!}y^{k}\bigg(\frac{x^{k}-1}{x-1}\bigg) \label{8} \\
&=\sum_{k=0}^{\infty}\frac{(1)_{k+1,\lambda}}{(k+1)!}y^{k+1}\sum_{n=0}^{k}x^{n}=\sum_{n=0}^{\infty}x^{n}\sum_{k=n+1}^{\infty}\frac{(1)_{k,\lambda}}{k!}y^{k}\nonumber \\
&=\sum_{n=0}^{\infty}x^{n}\bigg(e_{\lambda}(y)-1-\frac{(1)_{1,\lambda}}{1!}y-\frac{(1)_{2,\lambda}}{2!}y^{2}\cdots -\frac{(1)_{n,\lambda}}{n!}y^{n}\bigg)\nonumber.
\end{align}

Taking the limit as $x \rightarrow 1$ in \eqref{8}, we have 
\begin{align} 
&\sum_{n=0}^{\infty}\bigg(e_{\lambda}(y)-1-\frac{(1)_{1,\lambda}}{1!}y-\frac{(1)_{2,\lambda}}{2!}-\cdots-\frac{(1)_{n,\lambda}}{n!}y^{n}\bigg)\label{9} \\
&= \lim_{x \rightarrow 1}\sum_{k=1}^{\infty}\frac{(1)_{k,\lambda}}{k!}y^{k}\bigg(\frac{x^{k}-1}{x-1}\bigg)  =\sum_{k=1}^{\infty}\frac{(1)_{k,\lambda}}{k!}y^{k}k \nonumber \\
&=y\sum_{k=0}^{\infty}\frac{(1-\lambda)_{k,\lambda}}{k!}y^{k}
=ye_{\lambda}^{1-\lambda}(y)=\frac{y}{1+\lambda y}e_{\lambda}(y). \nonumber	
\end{align}
Therefore, by \eqref{8} and \eqref{9}, we obtain the following theorem. 
\begin{theorem}
The following identities hold true. 
\begin{align*}
&\frac{1}{x-1}\big(e_{\lambda}(xy)-e_{\lambda}(y)\big)\\
&=\sum_{n=0}^{\infty}\bigg(e_{\lambda}(y)-1-\frac{(1)_{1,\lambda}}{1!}y-\frac{(1)_{2,\lambda} }{2!}y^{2}-\cdots-\frac{(1)_{n,\lambda}}{n!}y^{n}\bigg)x^{n}, \\
&\frac{y}{1+\lambda y}e_{\lambda}(y)\\
&=\sum_{n=0}^{\infty}\bigg(e_{\lambda}(y)-1-\frac{(1)_{1,\lambda}}{1!}y-\frac{(1)_{2,\lambda}}{2!}y^{2}-\cdots-\frac{(1)_{n,\lambda}}{n!}y^{n}\bigg).
\end{align*}
\end{theorem}

The degenerate hyperbolic cosine is defined by 
\begin{equation*}
\cosh_{\lambda}(x)=\frac{e_{\lambda}(-x)+e_{\lambda}(x)}{2}.	
\end{equation*}
Note that $\lim_{\lambda\rightarrow 0}\cosh_{\lambda}(x)=\cosh(x)$.
The next corollary is immediate from Theorem 2.1.

\begin{corollary}
The following identities hold true. 
\begin{align*}
&\sum_{n=1}^{\infty}\bigg(e_{\lambda}(1)-1-\frac{(1)_{1,\lambda}}{1!}-\frac{(1)_{2,\lambda}}{2!}-\cdots-\frac{(1)_{n,\lambda}}{n!}\bigg)x^{n}=\frac{e_{\lambda}(x)-xe_{\lambda}(1)}{x-1}+1, \\
&\sum_{n=1}^{\infty}\bigg(e_{\lambda}(1)-1-\frac{(1)_{1,\lambda}}{1!}-\frac{(1)_{2,\lambda}}{2!}-\cdots-\frac{(1)_{n,\lambda}}{n!}\bigg)=1-\frac{\lambda}{1+\lambda}e_{\lambda}(1),
\end{align*}
and 
\begin{displaymath}
\sum_{n=1}^{\infty}\bigg(e_{\lambda}(1)-1-\frac{(1)_{1,\lambda}}{1!}-\frac{(1)_{2,\lambda}}{2!}-\cdots-\frac{(1)_{n,\lambda}}{n!}\bigg)(-1)^{n}=1-\cosh_{\lambda}(1). 
\end{displaymath}
\end{corollary}

From \eqref{7}, we note that 
\begin{align}
&\sum_{n=0}^{\infty}\binom{n}{p}\bigg(e_{\lambda}(y)-1-\frac{(1)_{1,\lambda}}{1!}y-\frac{(1)_{2,\lambda}}{2!}y^{2}-\cdots-\frac{(1)_{n,\lambda}}{n!}y^{n}\bigg) \label{10} \\
&=\sum_{n=p}^{\infty}\binom{n}{p}\sum_{k=n+1}^{\infty}\frac{(1)_{k,\lambda}}{k!}y^{k}=\sum_{k=p+1}^{\infty}\frac{(1)_{k,\lambda}}{k!}y^{k}\sum_{n=p}^{k-1}\binom{n}{p}\nonumber \\
&=\sum_{k=p+1}^{\infty}\frac{(1)_{k,\lambda}}{k!}y^{k}\sum_{n=p}^{k-1}\bigg(\binom{n+1}{p+1}-\binom{n}{p+1}\bigg)=\sum_{k=p+1}^{\infty}\frac{(1)_{k,\lambda}}{k!}y^{k}\binom{k}{p+1} \nonumber \\
&=\sum_{k=0}^{\infty}\frac{(1)_{k+p+1,\lambda}}{(k+p+1)!}y^{k+p+1}\binom{k+p+1}{p+1}=\frac{y^{p+1}(1)_{p+1,\lambda}}{(p+1)!}\sum_{k=0}^{\infty}\frac{(1-(p+1)\lambda)_{k,\lambda}}{k!}y^{k}\nonumber\\
&=\frac{y^{p+1}}{(p+1)!}(1)_{p+1,\lambda}e_{\lambda}^{1-(p+1)\lambda}(y)=\frac{y^{p+1}}{(p+1)!}(1)_{p+1,\lambda}(1+\lambda y)^{-(p+1)}e_{\lambda}(y).  \nonumber 
\end{align}
Therefore, by \eqref{10}, we obtain the following theorem. 
\begin{theorem}
For $p\ge 0$, we have 
\begin{align*}
\sum_{n=0}^{\infty}&\binom{n}{p}\bigg(e_{\lambda}(y)-1-\frac{(1)_{1,\lambda}}{1!}y-\frac{(1)_{2,\lambda}}{2!}y^{2}-\cdots-\frac{(1)_{n,\lambda}}{n!}y^{n}\bigg) \nonumber \\
&=\frac{y^{p+1}}{(p+1)!}(1)_{p+1,\lambda}(1+\lambda y)^{-(p+1)}e_{\lambda}(y). 
\end{align*}
Especially, for $y=1$, we obtain
\begin{align*}
\sum_{n=0}^{\infty}&\binom{n}{p}\bigg(e_{\lambda}(1)-1-\frac{(1)_{1,\lambda}}{1!}-\frac{(1)_{2,\lambda}}{2!}-\cdots-\frac{(1)_{n,\lambda}}{n!}\bigg) \\ &=\frac{(1)_{p+1,\lambda}}{(p+1)!}(1+\lambda)^{-(p+1)}e_{\lambda}(1). 
\end{align*}
\end{theorem}
From Theorem 2.3, we note that 
\begin{align}
\sum_{n=0}^{\infty}&(n)_{p}\bigg(e_{\lambda}(y)-1-\frac{(1)_{1,\lambda}}{1!}y-\frac{(1)_{2,\lambda}}{2!}y^{2}-\cdots-\frac{(1)_{n,\lambda}}{n!}y^{n}\bigg) \label{11}\\
&=\frac{y^{p+1}}{p+1}(1)_{p+1,\lambda}(1+\lambda y)^{-(p+1)}e_{\lambda}(y). \nonumber
\end{align}
By \eqref{3} and \eqref{11}, we get 
\begin{align}
&\sum_{n=0}^{\infty}(n)_{p,\lambda}\bigg(e_{\lambda}(y)-1-\frac{(1)_{1,\lambda}}{1!}y-\frac{(1)_{2,\lambda}}{2!}y^{2}-\cdots-\frac{(1)_{n,\lambda}}{n!}y^{n}\bigg) \label{12}\\
&=\sum_{n=0}^{\infty}\sum_{k=0}^{p}S_{2,\lambda}(p,k)(n)_{k}\bigg(e_{\lambda}(y)-1-\frac{(1)_{1,\lambda}}{1!}y-\frac{(1)_{2,\lambda}}{2!}y^{2}-\cdots-\frac{(1)_{n,\lambda}}{n!}y^{n}\bigg) \nonumber\\
&=\sum_{k=0}^{p}S_{2,\lambda}(p,k)\sum_{n=0}^{\infty}(n)_{k}\bigg(e_{\lambda}(y)-1-\frac{(1)_{1,\lambda}}{1!}y-\frac{(1)_{2,\lambda}}{2!}y^{2}-\cdots-\frac{(1)_{n,\lambda}}{n!}y^{n}\bigg) \nonumber\\
&=\sum_{k=0}^{p}S_{2,\lambda}(p,k)\frac{y^{k+1}}{k+1}(1)_{k+1,\lambda}(1+\lambda y)^{-(k+1)}e_{\lambda}(y). \nonumber
\end{align}
Therefore, by \eqref{12}, we obtain the following theorem. 
\begin{theorem}
For $p\ge 0$, we have 
\begin{align*}
&\sum_{n=0}^{\infty}(n)_{p,\lambda}\bigg(e_{\lambda}(y)-1-\frac{(1)_{1,\lambda}}{1!}y-\frac{(1)_{2,\lambda}}{2!}y^{2}-\cdots-\frac{(1)_{n,\lambda}}{n!}y^{n}\bigg) \\
&=\sum_{k=0}^{p}S_{2,\lambda}(p,k)\frac{y^{k+1}}{k+1}(1)_{k+1,\lambda}(1+\lambda y)^{-(k+1)}e_{\lambda}(y).
\end{align*}
In particular, for $y=1$, we get 
\begin{align*}
&\sum_{n=0}^{\infty}(n)_{p,\lambda}\bigg(e_{\lambda}(1)-1-\frac{(1)_{n,\lambda}}{1!}-\frac{(1)_{2,\lambda}}{2!}-\cdots-\frac{(1)_{n,\lambda}}{n!}\bigg) \\
&=\sum_{k=0}^{p}S_{2,\lambda}(p,k)\frac{(1)_{k+1,\lambda}}{k+1}(1+\lambda)^{-(k+1)}e_{\lambda}(1). 
\end{align*}
\end{theorem}
From \eqref{4}, we note that 
\begin{align}
&\sum_{n=0}^{\infty}S_{2,\lambda}(n,k)\frac{t^{n}}{n!}=\sum_{n=k}^{\infty}S_{2,\lambda}(n,k)\frac{t^{n}}{n!}=\frac{1}{k!}\big(e_{\lambda}(t)-1\big)^{k}\label{13} \\
&=\frac{1}{k!}\sum_{j=0}^{k}\binom{k}{j}(-1)^{k-j}e_{\lambda}^{j}(t)=\frac{1}{k!}\sum_{j=0}^{k}\binom{k}{j}(-1)^{k-j}\sum_{n=0}^{\infty}\frac{(j)_{n,\lambda}}{n!}t^{n} \nonumber\\
&=\sum_{n=0}^{\infty}\bigg(\frac{1}{k!}\sum_{j=0}^{k}\binom{k}{j}(-1)^{k-j}(j)_{n,\lambda}\bigg)\frac{t^{n}}{n!}.\nonumber	
\end{align}
Comparing the coefficients on both sides of \eqref{13}, we obtain 
\begin{equation}
S_{2,\lambda}(n,k)=\frac{1}{k!}\sum_{j=0}^{k}\binom{k}{j}(-1)^{k-j}(j)_{n,\lambda},\quad (n,k\ge 0).\label{14}
\end{equation}
Taking the limit as $\lambda \rightarrow 0$ in \eqref{14}, we have
\begin{equation}
S_{2}(n,k)=\frac{1}{k!}\sum_{j=0}^{k}\binom{k}{j}(-1)^{k-j}j^{n},\quad (n,k\ge 0).\label{15}
\end{equation}
By using \eqref{15}, we derive the following: 
\begin{align}
&\frac{1}{k!}\sum_{j=0}^{k}\binom{k}{j}(-1)^{k-j}\frac{e_{\lambda}(jy)-e_{\lambda}(y)}{j-1}\label{16}\\
&=\frac{1}{k!}\sum_{j=0}^{k}\binom{k}{j}(-1)^{k-j}\sum_{n=1}^{\infty}\frac{(1)_{n,\lambda}}{n!}y^{n}\bigg(\frac{j^{n}-1}{j-1}\bigg) \nonumber\\
&=\frac{1}{k!}\sum_{j=0}^{k}\binom{k}{j}(-1)^{k-j}\sum_{n=1}^{\infty}\frac{(1)_{n,\lambda}}{n!}y^{n}\sum_{l=0}^{n-1}j^{l} \nonumber \\
&=\sum_{n=1}^{\infty}\frac{(1)_{n,\lambda}}{n!}y^{n}\sum_{l=0}^{n-1}\frac{1}{k!}\sum_{j=0}^{k}\binom{k}{j}(-1)^{k-j}j^{l} \nonumber\\
&=\sum_{n=1}^{\infty}\frac{(1)_{n,\lambda}}{n!}y^{n}\sum_{l=0}^{n-1}S_{2}(l,k)=\sum_{l=0}^{\infty}S_{2}(l,k)\sum_{n=l+1}^{\infty}\frac{(1)_{n,\lambda}}{n!}y^{n}.\nonumber
\end{align}
Therefore, by \eqref{16}, we obtain the following theorem. 
\begin{theorem}
For $k\ge 0$, we have 
\begin{align*}
\sum_{n=0}^{\infty}&S_{2}(n,k)\bigg(e_{\lambda}(y)-1-\frac{(1)_{1,\lambda}}{1!}y-\frac{(1)_{2,\lambda}}{2!}y^{2}-\cdots-\frac{(1)_{n,\lambda}}{n!}y^{n}\bigg) \\
&=\frac{1}{k!}\sum_{j=0}^{k}\binom{k}{j}(-1)^{k-j}\frac{e_{\lambda}(jy)-e_{\lambda}(y)}{j-1}.
\end{align*}
In particular, for $y=1$, we get 
\begin{align*}
\sum_{n=0}^{\infty}&S_{2,\lambda}(n,k)\bigg(e_{\lambda}(1)-1\frac{(1)_{1,\lambda}}{1!}-\frac{(1)_{2,\lambda}}{2!}-\cdots-\frac{(1){n,\lambda}}{n!}\bigg)\\
&=\frac{1}{k!}\sum_{j=0}^{k}\binom{k}{j}(-1)^{k-j}\frac{e_{\lambda}(j)-e_{\lambda}(1)}{j-1}.
\end{align*}
\end{theorem}
\begin{remark}
We may naturally consider the following problem. \par 
\noindent For any $k\ge 0$, find the value of 
\begin{displaymath}
\sum_{n=0}^{\infty}S_{2,\lambda}(n,k)\bigg(e_{\lambda}(y)-1-\frac{(1)_{1,\lambda}}{1!}y-\frac{(1)_{2,\lambda}}{2!}y^{2}-\cdots-\frac{(1)_{n,\lambda}}{n!}y^{n}\bigg).
\end{displaymath}
\end{remark}
\begin{remark} 
Much work has been done as to degenerate and truncated theories. These theories have some applications to mathematics, engineering and physics. Researchers interested in these may refer to [1-22]. 
\end{remark}

\section{Conclusion}
In this note, we studied infinite series whose terms involve the truncated degenerate exponentials together with binomial coefficients, the generalized falling factorials and the degenerate Stirling numbers of the second kind and determined either their values or some other expressions of them as finite sums. \par
In recent years, we have witnessed that study of degenerate versions yielded many fascinating and fruitful results. We would like to continue to study degenerate versions of many special numbers and polynomials and to find some applications of them to physics, science and engineering.

\end{document}